\documentclass[11pt]{article}

\usepackage{amsmath,amsthm,amsfonts,epsf}
\theoremstyle{plain}
\newtheorem{theorem}{Theorem}[section]
\numberwithin{equation}{section}

\newcommand{\qbin}[2]{\genfrac{[}{]}{0pt}{}{#1}{#2}}
\newcommand{\qbins}[2]{{\textstyle\genfrac{[}{]}{0pt}{}{#1}{#2}}}
\newcommand{\Z}{\mathbb{Z}}
\newcommand{\I}{\mathcal{I}}
\newcommand{\T}{\mathcal{T}}
\newcommand{\A}{\text{A}}
\newcommand{\D}{\text{D}}
\newcommand{\E}{\text{E}}
\newcommand{\G}{\mathfrak{g}}

\begin{document}

\title{Refined $q$-trinomial coefficients and character identities}
 
\author{\large
S.~Ole Warnaar\thanks{
Instituut voor Theoretische Fysica,
Universiteit van Amsterdam,
Valckenierstraat 65,
1018 XE Amsterdam, The Netherlands;
e-mail: {\tt warnaar@wins.uva.nl}}}

\date{\small
Dedicated to Rodney Baxter on the occasion of his sixtiest birthday} 
\maketitle

\begin{abstract}
A refinement of the $q$-trinomial coefficients is introduced,
which has a very powerful iterative property.
This ``$\T$-invariance'' is applied to derive new Virasoro character
identities related to the exceptional simply-laced Lie algebras
$\E_6,\E_7$ and $\E_8$. 
\end{abstract}

\section{Introduction}
\subsection{Rodney Baxter}
Rodney Baxter is justly famous for his many beautiful discoveries in
mathematics and physics.
The 8-vertex model, Yang--Baxter equation, corner transfer matrix and
hard-hexagon model are among his most envied mathematical trophies.
This paper deals with a less-well-known discovery
of Rodney Baxter (made together with George Andrews),
that of the $q$-trinomial coefficients~\cite{AB87}.
My main aim will be to (for once) prove Baxter (and Andrews) wrong,
and show that the statement
\begin{quote}
``The literature is sparse on trinomial coefficients perhaps because they
lack both depth and elegance.'',
\end{quote}
made in the introduction of \cite{AB87}, is not at all justified.

To have any chance of succeeding, I have omitted all proofs
in this paper (which \textit{are} lacking elegance indeed!).
These will be given in a forthcoming longer paper on the same 
topic.

\subsection{$q$-Trinomial coefficients}
In their joined work on a generalization of the hard-hexagon model,
Andrews and Baxter \cite{AB87} were led to consider $q$-deformations
of the numbers appearing in the following generalized Pascal triangle:
\begin{equation*}
\begin{array}{ccccccccc}
&&&&1&&&& \\
&&&1&1&1&&& \\
&&1&2&3&2&1&& \\
&1&3&6&7&6&3&1& \\
.&.&.&.&.&.&.&.&.
\end{array}
\end{equation*}
The generating function for the numbers appearing in the $(L+1)$th row
is $(1+x+x^2)^L$, so that an explicit expression for the trinomial
coefficients can be found by double application of Newton's binomial
expansion. Explicitly,
\begin{equation*}
(1+x+x^2)^L=\sum_{a=-L}^L\binom{L}{a}_2 x^{a+L},
\end{equation*}
with
\begin{equation*}
\binom{L}{a}_2=\sum_{k\geq 0} \binom{L}{k}\binom{L-k}{k+a}.
\end{equation*}
(The effective range of summation is from $\max\{0,-a\}$ to 
$\min\{L,\lfloor(L-a)/2\rfloor\}$ so that one indeed finds
a nonzero number for $|a|\leq L$ only.)

Andrews and Baxter introduced several $q$-analogues of the trinomial
coefficients. Here we shall restrict ourselves to the simplest two given 
by~\cite[Eq. (2.7); $B=A$]{AB87}
\begin{equation*}
\qbin{L;q}{a}_2=\qbin{L}{a}_2
=\sum_{k\geq 0} q^{k(k+a)}\qbin{L}{k}\qbin{L-k}{k+a}
\end{equation*}
and
\begin{equation}\label{qtT}
T(L,a;q)=T(L,a)=
q^{\frac{1}{2}(L-a)(L+a)}\qbin{L;q^{-1}}{a}_2.
\end{equation}
(This is $T_0(L,a;q^{1/2})$ of \cite{AB87}.)
Here
\begin{equation*}
\qbin{n}{a}=
\begin{cases}\displaystyle \frac{(q)_{n}}{(q)_a(q)_{n-a}} &
\text{for $0\leq a \leq n$} \\[3mm]
0 & \text{otherwise,} \end{cases}
\end{equation*}
is a $q$-binomial coefficient or Gaussian polynomial, with
$(a;q)_n=(a)_n=\prod_{j=0}^{n-1} (1-aq^j)$ for $n\geq 1$ and
$(a;q)_0=(a)_0=1$. A convenient explicit expression for $T(L,a)$
is given by~\cite[Eq. (2.60)]{AB87}
\begin{equation*}
T(L,a)=\sum_{\substack{n=0 \\n+a+L \text{ even}}}^{L-|a|}
\frac{q^{\frac{1}{2}n^2}(q)_L}{(q)_{\frac{L-a-n}{2}}
(q)_{\frac{L+a-n}{2}}(q)_n}\, .
\end{equation*}

Some useful properties of the $q$-trinomial coefficients are
the symmetries $\qbins{L}{a}_2=\qbins{L}{-a}_2$ and $T(L,a)=T(L,-a)$,
and the large $L$ limits
\begin{equation}\label{tlim}
\lim_{L\to\infty}\qbin{L}{a}_2=\frac{1}{(q)_{\infty}}
\end{equation}
and
\begin{equation}\label{Tlim}
\lim_{\substack{L\to\infty \\ L+a+\sigma \text{ even}}}T(L,a)=
\sum_{\substack{n=0 \\ n+\sigma \text{ even}}}^{\infty}
\frac{q^{\frac{1}{2}n^2}}{(q)_n}=c_{\sigma}(q).
\end{equation}
Here $c_0$ and $c_1$ are (normalized) level-1 string functions of
$\A_1^{(1)}$, which admit the alternative representations
\begin{align}\label{sf}
c_{\sigma}(q)&=
\frac{(-q^{1/2};q)_{\infty}+(-1)^{\sigma}(q^{1/2};q)_{\infty}}
{2(q;q)_{\infty}} \\
&=\frac{q^{\frac{1}{2}\sigma}}{(q;q)_{\infty}
(q^{3-2\sigma},q^4,q^{5+2\sigma};q^8)_{\infty}
(q^{2+4\sigma},q^{14-4\sigma};q^{16})_{\infty}}  \notag
\end{align}
with the convention that $(a_1,\dots,a_k;q)_n=(a_1;q)_n\dots(a_k;q)_n$.

\section{A refinement of the $q$-trinomial coefficients}

For integers $L,M,a$ and $b$ we define the polynomial
\begin{multline*}
\T(L,M,a,b;q)=\T(L,M,a,b) \\
=\sum_{\substack{n=0 \\n+a+L \text{ even}}}^{\min\{L-|a|,M\}}
q^{\frac{1}{2}n^2}
\qbins{M}{n}
\qbins{M+b+(L-a-n)/2}{M+b}
\qbins{M-b+(L+a-n)/2}{M-b}.
\end{multline*}
Some trivial properties of $\T$ are
\begin{equation*}
\T(L,M,a,b)=0 \quad \text{if $|a|>L$ or $|b|>M$}
\end{equation*}
(the if is not an iff),
the symmetry
\begin{equation*}
\T(L,M,a,b)=\T(L,M,-a,-b),
\end{equation*}
the duality
\begin{equation}\label{mTdual}
\T(L,M,a,b;1/q)=q^{ab-ML}\T(L,M,a,b;q)
\end{equation}
and the limit
\begin{equation}\label{mTlim}
\lim_{M\to\infty}\T(L,M,a,b)=\frac{T(L,a)}{(q)_L}.
\end{equation}
What is perhaps less evident is that $\T$ can be viewed
as a refinement of both types of $q$-trinomial coefficients
in the following sense:
\begin{equation}\label{mTtoT}
\sum_{i=|b|}^{L-|a-b|}q^{\frac{1}{2}(i^2-b^2)}\T(L-i,i,a-b,b)=T(L,a)
\end{equation}
and
\begin{equation}\label{mTtot}
\sum_{i=|b|}^{L-|a-b|}q^{\frac{1}{2}(i^2-b^2)} \T(i,L-i,b,a-b)=\qbin{L}{a}_2.
\end{equation}
Here it is assumed that $a\geq b\geq 0$ or $a\leq b\leq 0$ in both formulas.
We note that the second equation follows from the first
by application of \eqref{qtT} and \eqref{mTdual}.
Equation \eqref{mTtoT} results after taking $M\to\infty$ in 
Theorem~\ref{mainthm} of the next section.

As an example of \eqref{mTtoT} let us calculate $T(4,2)$ in three different
ways. When $b=0$ in \eqref{mTtoT} we get
\begin{align*}
T(4,2)&=\T(4,0,2,0)+q^{1/2}\T(3,1,2,0)+q^2\T(2,2,2,0) \\
&=1+q(1+q+q^2)+q^2(1+q+2q^2+q^3+q^4),
\end{align*}
when $b=1$,
\begin{align*}
T(4,2)&=\T(3,1,1,1)+q^{3/2}\T(2,2,1,1)+q^4\T(1,3,1,1) \\
&=1+q+q^2+q^2(1+q)^2+q^4(1+q+q^2)
\end{align*}
and, finally, when $b=2$,
\begin{align*}
T(4,2)&=\T(2,2,0,2)+q^{5/2}\T(1,3,0,2)+q^6\T(0,4,0,2) \\
&=1+q+2q^2+q^3+q^4+q^3(1+q+q^2)+q^6.
\end{align*}
Simplifying each of these three expressions correctly yields
$T(4,2)=1+q+2q^2+2q^3+2q^4+q^5+q^6$.

To conclude this section we remark that \eqref{mTtot} is a bounded analogue
of the following summation \cite[Eq. (4.3); $-q^{-n}\to\infty$]{AB98},
\cite[Eq. (2.10)]{BMP98}
\begin{equation}\label{ABP}
\sum_{i=0}^{\infty}q^{\frac{1}{2}i^2}\frac{T(i,b)}{(q)_i}=
\frac{q^{\frac{1}{2}b^2}}{(q)_{\infty}}
\end{equation}
as can be seen by taking $L$ to infinity in \eqref{mTtot}
using \eqref{tlim} and \eqref{mTlim}.
Hence \eqref{mTtot} should be compared with~\cite[Eq. (10)]{W99}
\begin{equation}\label{con}
\sum_{i=0}^{\infty}q^{\frac{1}{2}i^2}\qbin{L}{i}T(i,b)=
q^{\frac{1}{2}b^2}\qbin{2L}{L-b},
\end{equation}
which also yields \eqref{ABP} in the large $L$ limit.

\section{$\T$-invariance}
The important question to be addressed is whether the refined $q$-trinomial
$\T$ is at all relevant. The answer to this is a clear ``yes''.
Not only did we find that almost any result for $q$-trinomial coefficients
has an analogue for the polynomials $\T$ (\eqref{con} appears to be an 
exception), but, thanks to the following theorem, $\T$ has perhaps even 
more depth than the $q$-trinomials.
\begin{theorem}\label{mainthm}
For $L,M,a,b$ integers such that $a,b\geq 0$ or $a,b\leq 0$ there holds
\begin{equation}\label{Tinv}
\sum_{i=|b|}^{\min\{L-|a|,M\}} q^{\frac{1}{2}i^2}
\qbin{L+M-i}{L}\T(L-i,i,a,b)
= q^{\frac{1}{2}b^2}\T(L,M,a+b,b).
\end{equation}
\end{theorem}
This is a very powerful summation formula that allows one to iterate
identities involving $\T$, thereby generating an infinite chain
of $\T$-identities.
By the limit \eqref{mTlim} this then produces an infinite chain of
$q$-trinomial identities, and hence (by \eqref{qtT}--\eqref{Tlim})
of $q$-series identities.

\section{Three exceptional examples}
Taking simple $\T$-identities such as
$$\sum_{j\in\Z}q^{j(j+1)}\{\T(L,M,2j,j)-\T(L,M,2j+2,j)\}=
\delta_{L,0}\delta_{M,0}$$
or
$$\sum_{j\in\Z}(-1)^j q^{j(j+1)/2}\{\T(L,M,j,j)-\T(L,M,j+1,j)\}=
\delta_{M,0}$$
as input, and using the $\T$-invariance of Theorem~\ref{mainthm} 
to iterate these,
we have proved large classes of identities for doubly bounded analogues of
Virasoro characters. Many of the limiting character identities are known and
many are new. 
In this paper we shall, however, not prove a single identity 
using \eqref{Tinv}.
Instead, we shall only try to demonstrate the power of 
the theorem to generate new identities. As input we take three
identities for which, at present, we have not a clue to a proof.
However, accepting these initial conjectures, six beautiful series
of identities follow, which would have been almost impossible
to conceive without Theorem~\ref{mainthm}.

\subsection{Some preliminaries}

\begin{figure}[hbt]
\centerline{\epsfxsize=10 cm \epsffile{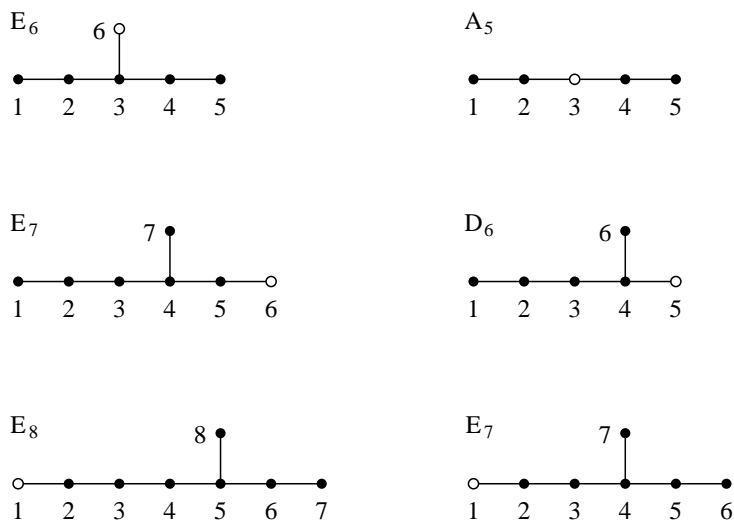}}
\caption{\small Dynkin diagrams of the Lie algebras $\E_n$, $\A_5$ and $\D_6$
with labelling of vertices as used in the text.
Removing the marked vertex (and its corresponding edge) in a diagram
of the first column yields the diagram to its right. Conversely, adding
a vertex plus edge to the marked vertices in the second column yields
the graphs of the first column.}
\end{figure}

Let $\G$  be any of the simply laced Lie algebras
whose Dynkin diagram is shown in Figure~1.
Given the Dynkin diagram of $\G$ together with its labelling of vertices,
we define a corresponding incidence matrix $\I_{\G}$
with entries
\begin{equation*}
(\I_{\G})_{i,j}=\begin{cases}
1 & \text{if vertices $i$ and $j$ are connected by an edge} \\
0 & \text{otherwise,}
\end{cases}
\end{equation*}
where $i,j=1,\dots,r_{\G}$ with $r_{\G}$ the rank of $\G$ (i.e., the
number of vertices of the diagram).
For any $\G$ we define an $(m,n)$-system
as the set of $r_{\G}$ linear coupled equations
\begin{equation}\label{mn}
m+n=\frac{1}{2}(\I_{\G}m+Ne_i).
\end{equation}
Here $N$ is a nonnegative integer, $m,n,e_i$ are vectors in
$\Z_+^{r_{\G}}$, with $e_i$ the unit vector ($(e_i)_j=\delta_{i,j}$)
associated with the $i$th vertex of the Dynkin diagram of $\G$.
Only those labels $i$ will occur that correspond to the marked vertices
(drawn as open circles) in Figure~1.
This fixes $i$ for all $\G$ other than $\E_7$.

For given $N$ and $i$, $m$ determines $n$ and vice versa. 
If $C_{\G}$ is the Cartan matrix of $\G$,
i.e., $C_{\G}=2I-\I_{\G}$, we find explicitly that
\begin{equation*}
n=\frac{1}{2}(N e_i-C_{\G}m) \quad \text{and} \quad 
m=C^{-1}_{\G}(N e_i-2n).
\end{equation*}
Note though that not all $m~(n)$ with integer entries
will also yield an $n~(m)$ with integer entries.

As an example let $\G=\E_7$, $N=6$ and $i=1$.
Then the only admissible solutions to \eqref{mn} are
\begin{align}\label{mnsolneven}
m&=5e_1+4e_2+3e_3+2e_4+e_7, & n&=e_5 \notag \\
m&=3e_1+4e_2+5e_3+6e_4+4e_5+2e_6+3e_7, & n&=2e_1 \notag \\
m&=5e_1+4e_2+5e_3+6e_4+4e_5+2e_6+3e_7, & n&=e_2 \\
m&=7e_1+8e_2+9e_3+10e_4+6e_5+2e_6+5e_7, & n&=e_6 \notag \\
m&=9e_1+12e_2+15e_3+18e_4+12e_5+6e_6+9e_7, & n&=0 \notag 
\end{align}
and
\begin{align}\label{mnsolnodd}
m&=0, & n&=3e_1 \notag \\
m&=2e_1, & n&=e_1+e_2 \notag \\
m&=4e_1+2e_2, & n&=e_3 \notag \\
m&=4e_1+4e_2+4e_3+4e_4+2e_5+2e_7, & n&=e_1+e_6 \\
m&=6e_1+6e_2+6e_3+6e_4+4e_5+2e_6+2e_7, & n&=e_7 \notag \\
m&=6e_1+8e_2+10e_3+12e_4+8e_5+4e_6+6e_7, & n&=e_1, \notag 
\end{align}
where the first (second) set of solutions meets the criterion that
$n_1+n_3+n_7$ is even (odd).

Finally we need polynomials associated with the algebras $\G$
depicted in the second column of Figure~1 as follows.
Let $p=3,5,1$ for $\G=\A_5,\D_6,\E_7$, respectively,
so that $p$ corresponds to the marked vertex of $\G$.
For $M$ a nonnegative integer, $\sigma=0,1$ and $(m,n)$-system 
\begin{equation}\label{mnG}
m+n=\frac{1}{2}(\I_{\G}m+2Me_p)
\end{equation}
we define
\begin{equation*}
F^{\A_5}_{M;\sigma}(q)=
\sum_{\substack{n\in\Z^5_+ \\ n_1+n_4 \equiv n_2+n_5 \pmod{3}
\\n_1+n_3+n_5+\sigma \text{ even}}}
q^{nC^{-1}_{\A{5}}n}\qbin{m+n}{n}
\end{equation*}
and
\begin{equation*}
F^{\D_6}_{M;\sigma}(q)=
\sum_{\substack{n\in\Z^6_+ \\ n_1+n_3+n_6 \text{ even}
\\n_1+n_3+n_5+\sigma \text{ even}}}
q^{nC^{-1}_{\D{6}}n}\qbin{m+n}{n}
\end{equation*}
and
\begin{equation*}
F^{\E_7}_{M;\sigma}(q)=
\sum_{\substack{n\in\Z^7_+ \\ n_1+n_3+n_7+\sigma \text{ even}}}
q^{nC^{-1}_{\E{7}}n}\qbin{m+n}{n}.
\end{equation*}
Here we have used the abbreviations
$nC^{-1}_{\G}n=\sum_{i,j=1}^{r_{\G}}(C^{-1}_{\G})_{i,j} n_i n_j$
and
$\qbins{m+n}{n}=\prod_{j=1}^{r_{\G}}\qbins{m_j+n_j}{n_j}$
for $m,n\in\Z^{r_{\G}}$.
Similarly we will write $(q)_n=\prod_{j=1}^{r_{\G}}(q)_{n_j}.$

As example consider again $\G=\E_7$ and choose $M=3$. 
Then the only contributing terms to $F^{\E_7}_{3;0}(q)$ and 
$F^{\E_7}_{3;1}(q)$ correspond to the solutions of \eqref{mn}
listed in \eqref{mnsolneven} and \eqref{mnsolnodd}, respectively. Hence
$$F^{\E_7}_{3;0}(q)=q^6+q^6\qbins{5}{2}+q^4\qbins{5}{1}+q^2\qbins{3}{1}+1$$
and
$$F^{\E_7}_{3;1}(q)=
q^{27/2}+q^{19/2}\qbins{3}{1}+q^{15/2}+q^{11/2}\qbins{5}{1}
+q^{7/2}\qbins{3}{1}+q^{3/2}\qbins{7}{1}.$$

\subsection{An ($\E_7$,$\E_8$) series}
Our first conjecture is the following polynomial identity 
involving $\E_7$:
\begin{multline}\label{conj1}
\sum_{j\in\Z}\Bigl\{q^{\frac{1}{2}j(15j+2)}\T(L,M,3j,5j)
-q^{\frac{1}{2}(3j+1)(5j+1)}\T(L,M,3j+1,5j+1)\Bigr\} \\
=\sum_{\substack{n\in\Z^7_+ \\ n_1+n_3+n_7+L \text{ even}}}
q^{nC^{-1}_{\E_7}n}\qbin{\frac{1}{2}(L+M+m_1)}{2M}\qbin{m+n}{n},
\end{multline}
with $(m,n)$-system \eqref{mnG} and $\G=\E_7$ (so that $p=1$).
The restriction on the sum in the right side guaranties that, 
given $n\in\Z_+^7$, $L+M+m_1$ is even. 
Using Theorem~\ref{mainthm} to iterate this conjecture we obtain an
infinite series of polynomial identities.
These identities are best expressed by turning $\E_7$ into $\E_8$
by the mechanism described in the caption of Figure~1.
Specifically, for $k\geq 1$,
\begin{multline}\label{flower}
\sum_{j\in\Z}\Bigl\{q^{\frac{1}{2}j(5(5k+3)j+2)}\T(L,M,(5k+3)j,5j) \\
-q^{\frac{1}{2}(5j+1)((5k+3)j+k+1)}\T(L,M,(5k+3)j+k+1,5j+1)\Bigr\} \\
=\sum_{r\in\Z^{k-1}_+}
\Bigl(\prod_{a=0}^{k-2} q^{\frac{1}{2}(r_a-r_{a+1})^2}
\qbins{r_{a-1}-r_a+r_{a+1}}{r_a}\Bigr)
\sum_{n\in\Z^8_+}
q^{\frac{1}{4}mC_{\E_8}m}\qbins{r_{k-2}-\frac{1}{2}m_1}{r_{k-1}}
\qbins{m+n}{n},
\end{multline}
with $r_0=L$, $r_{-1}=L+M$ and $(m,n)$-system given by
\begin{equation}\label{mnE8}
m+n=\frac{1}{2}(\I_{\E_8}m+r_{k-1}e_1).
\end{equation}

We now use this result to obtain $q$-series identities.
First we let $M$ tend to infinity, which, by \eqref{mTlim}, 
turns \eqref{flower} into an identity for $q$-trinomial coefficients.
Then we either send $L$ to infinity using \eqref{Tlim}, or we
first replace $q\to 1/q$ and then send $L$ to infinity using \eqref{tlim}.
Omitting the actual calculations (which sometimes require variable changes
to remove $L$-dependent terms in the exponent of $q$)
we find two families of $q$-series identities, one 
of $\E_7$ type and one of $\E_8$ type.
To present these identities in a neat form we recall the bosonic 
representation of the Virasoro characters~\cite{FF83,RC85}
\begin{equation}\label{Vir}
\chi_{r,s}^{(p,p')}(q)=\frac{q^{\frac{(p'r-ps)^2-1}{4pp'}}}{(q)_{\infty}}
\sum_{j=-\infty}^{\infty}\Bigl\{
q^{j(pp'j+p'r-ps)}-q^{(pj+r)(p'j+s)}\Bigr\},
\end{equation}
for coprime integers $2\leq p<p'$ and $r=1,\dots,p-1$, $s=1,\dots,p'-1$.
The somewhat unusual normalization of the Virasoro characters is chosen to
simplify subsequent equations.
Besides the Virasoro characters we also need the following (subset)
of the branching functions corresponding to the coset
$(\A_1^{(1)}\oplus\A_1^{(1)},
\A_1^{(1)})$ at levels $2p/(p'-p)-2,2$ and $2p/(p'-p)$~\cite{KW90}:
\begin{multline}
B_{r,s;\sigma}^{(p,p')}(q)=\frac{q^{\frac{(p'r-ps)^2-4}{8pp'}}}{(q)_{\infty}}
\sum_{j=-\infty}^{\infty}\Bigl\{
q^{j(pp'j+p'r-ps)}c_{pj+\frac{1}{2}(r-s)+\sigma}(q) \\
-q^{(pj+r)(p'j+s)}c_{pj+\frac{1}{2}(r+s)+\sigma}(q)\Bigr\},
\end{multline}
for integers $p,p',r,s$ in the same ranges as above, such that
$p'-p$ and $r-s$ are even, and $\gcd((p'-p)/2,p')=1$. The integer $\sigma$
takes the values $0$ or $1$ and
the $c_j$ are the string functions of equations \eqref{Tlim} and
\eqref{sf} with the obvious identification of $c_{2j+\sigma}$ with
$c_{\sigma}$.

The various large $L$ and $M$ limits now lead to the following 
list of identities.
\begin{itemize}
\item
Taking $L+\sigma$ even in \eqref{flower}, sending $M$ and $L$ to
infinity using that $c_{\sigma}(q)=\chi^{(3,4)}_{\sigma+1,1}(q)/
(q)_{\infty}$ yields for $k\geq 1$
\begin{equation*}
\sum_{n_1,\dots,n_k\geq 0}
\frac{q^{\frac{1}{2}(N_1^2+\cdots+N_k^2)}F^{\E_7}_{n_k;m_{\sigma}}(q)}
{(q)_{n_1}\cdots(q)_{n_{k-1}}(q)_{2n_k}}
=\begin{cases} \displaystyle
\chi_{\sigma+1,1}^{(3,4)}(q)\chi_{1,(k+1)/2}^{(5,(5k+3)/2)}(q)
& \text{$k$ odd} \\[4mm] \displaystyle
B_{1,k+1;\sigma}^{(5,5k+3)}(q) & \text{$k$ even,}
\end{cases}
\end{equation*}
where the following definitions have been employed:
\begin{equation}\label{Nm}
N_a=n_a+\cdots+n_k \quad \text{and} \quad
m_{\sigma}\equiv \sigma+\sum_{\substack{a=1 \\a \text{ odd}}}^k n_a
\pmod{2},
\end{equation}
for $m_{\sigma}\in\{0,1\}$.
There is a corresponding ``$k=0$'' identity obtained by
taking the same limit as above, but now in the initial
conjecture \eqref{conj1}.
Using the (from a $q$-series point of view nontrivial)
relation $B_{1,1;\sigma}^{(3,5)}=\chi_{2\sigma+1,1}^{(4,5)}$,
which follows from a symmetry of the A$_1^{(1)}$ branching functions,
we find the well-known $\E_7$ conjecture~\cite{KKMM93}
\begin{equation}\label{E7conj}
\chi_{2\sigma+1,1}^{(4,5)}(q)=
\sum_{\substack{n\in\Z^7_+ \\ n_1+n_3+n_7+\sigma \text{ even}}}
\frac{q^{nC^{-1}_{\E_7}n}}{(q)_n}.
\end{equation}

\item
If in \eqref{flower} we send $M$ to infinity, 
replace $q\to 1/q$ and then take the
limit of large $L$ we obtain for $k\geq 2$,
\begin{multline}\label{X}
\sum_{r\in\Z^{k-1}_+}
\sideset{}{'}\sum_{m\in\Z_+^8}
\frac{q^{\frac{1}{2}\sum_{a=1}^{k-1}(r_a-r_{a-1})^2}}
{(q)_{r_1}}
\Bigl(\prod_{a=2}^{k-1} 
\qbins{r_{a-1}-r_a+r_{a+1}}{r_a}\Bigr)
q^{\frac{1}{4}mC_{\E_8}m}
\qbins{m+n}{m} \\
=\begin{cases}
\chi^{((5k+3)/2,5k-2)}_{(k+1)/2,k}(q) & \text{$k$ odd} \\[3mm]
\chi^{(5k/2-1,5k+3)}_{k/2,k+1}(q) & \text{$k$ even,}
\end{cases}
\end{multline}
with $r_0=0$, $r_k=r_{k-1}-m_1/2$ and $(m,n)$-system \eqref{mnE8}.
The prime in the sum over $m$ denotes the restriction
$m_2\equiv m_4\equiv m_8\equiv r_{k-1}\pmod{2}$ with all
other $m_i$ being even.
When $k=1$ the resulting character formula takes a somewhat different form,
and one obtains the $\E_8$ identity \cite{KKMM93,WP94}
\begin{equation}\label{E8}
\chi^{(3,4)}_{1,1}(q)=
\sum_{n\in\Z_+^8}
\frac{q^{nC^{-1}_{\E_8}n}}{(q)_n}
=\frac{1}{(q^3,q^4,q^5;q^8)_{\infty}(q^2,q^{14};q^{16})_{\infty}}.
\end{equation}
\end{itemize}

\subsection{A ($\D_6$,$\E_7$) series}
In our second conjecture the role of $\E_7$ is taken over by $\D_6$,
\begin{multline}\label{conj2}
\sum_{j\in\Z}\Bigl\{q^{\frac{1}{2}j(24j+2)}\T(L,M,4j,6j)
-q^{\frac{1}{2}(4j+1)(6j+1)}\T(L,M,4j+1,6j+1)\Bigr\} \\
=\sum_{\substack{n\in\Z^6_+ \\ n_1+n_3+n_6 \text{ even}\\
n_1+n_3+n_5+L \text{ even}}}
q^{nC^{-1}_{\D_6}n}\qbin{\frac{1}{2}(L+M+m_5)}{2M}\qbin{m+n}{n},
\end{multline}
with $(m,n)$-system \eqref{mnG} where $\G=\D_6$ (and $p=5$).
The restriction that $n_1+n_3+n_5+L$ is even is necessary 
for $L+M+m_5$ to be even. The additional constraint on
$n_1+n_3+n_6$ is to avoid an extra
$$q\sum_{j\in\Z}\Bigl\{q^{\frac{1}{2}j(24j+14)}\T(L,M,4j+1,6j+2)
-q^{\frac{1}{2}(4j+3)(6j+1)}\T(L,M,4j+2,6j+3)\Bigr\}$$
on the left-hand side. Of course, this means we have actually two 
conjectures, but the case when $n_1+n_3+n_6$ is odd will not be pursued here.

Using Theorem~\ref{mainthm} to iterate the $\D_6$ conjecture
we obtain an infinite series of polynomial identities involving $\E_7$.
Specifically, for $k\geq 1$ there holds
\begin{multline}\label{flower2}
\sum_{j\in\Z}\Bigl\{q^{\frac{1}{2}j(6(6k+4)j+2)}\T(L,M,(6k+4)j,6j) \\
-q^{\frac{1}{2}(6j+1)((6k+4)j+k+1)}\T(L,M,(6k+4)j+k+1,6j+1)\Bigr\} \\
=\sum_{r\in\Z^{k-1}_+}
\Bigl(\prod_{a=0}^{k-2} q^{\frac{1}{2}(r_a-r_{a+1})^2}
\qbins{r_{a-1}-r_a+r_{a+1}}{r_a}\Bigr) \\
\times
\sum_{\substack{n\in\Z^7_+ \\ n_1+n_3+n_7 \text{ even}}}
q^{\frac{1}{4}mC_{\E_7}m}\qbins{r_{k-2}-\frac{1}{2}m_6}{r_{k-1}}
\qbins{m+n}{n},
\end{multline}
where $r_0=L$, $r_{-1}=L+M$ and 
\begin{equation}\label{mnE7}
m+n=\frac{1}{2}(\I_{\E_7}m+r_{k-1}e_6).
\end{equation}
As before we consider the various large $L$ and $M$ limits.
\begin{itemize}
\item
Taking $L+\sigma$ even in \eqref{flower2} and sending $M$ and $L$ to
infinity yields for $k\geq 1$ that
\begin{equation*}
\sum_{n_1,\dots,n_k\geq 0}
\frac{q^{\frac{1}{2}(N_1^2+\cdots+N_k^2)}F^{\D_6}_{n_k;m_{\sigma}}(q)}
{(q)_{n_1}\cdots(q)_{n_{k-1}}(q)_{2n_k}}
=\begin{cases} \displaystyle
\chi_{\sigma+1,1}^{(3,4)}(q)\chi_{1,(k+1)/2}^{(6,3k+2)}(q)
& \text{$k$ odd} \\[4mm]
\displaystyle B_{1,k+1;\sigma+k/2}^{(6,6k+4)}(q) & \text{$k$ even,}
\end{cases}
\end{equation*}
with the notation of equation \eqref{Nm} and with 
the identification of $B_{r,s;\sigma+2j}^{(p,p')}$
with $B_{r,s;\sigma}^{(p,p')}$.

The $k=0$ case, corresponding to the above limit taken in \eqref{conj2}
yields
\begin{equation*}
B_{1,1;\sigma}^{(4,6)}=
\sum_{\substack{n\in\Z^6_+ \\ n_1+n_3+n_6 \text{ even}\\
n_1+n_3+n_5 \equiv \sigma \text{ even}}}
\frac{q^{nC^{-1}_{\D_6}n}}{(q)_{n}}.
\end{equation*}
Although we were unable to prove this, it appears that the above branching
function admits the following simplification
\begin{equation*}
B_{1,1;\sigma}^{(4,6)}(q)=
\begin{cases}\displaystyle
\frac{1}{(q)_{\infty}}\Bigl(
\sum_{j=0}^{\infty}(-q)^{j^2}+\sum_{j=1}^{\infty}q^{6j^2}\Bigr)
&\text{$\sigma=0$} \\[5mm]
\displaystyle 
\frac{q^{3/2}}{(q)_{\infty}} \sum_{j=0}^{\infty}q^{6j(j+1)}
=\frac{q^{3/2}(q^{24};q^{24})_{\infty}}{(q^{12};q^{24})_{\infty}(q;q)_{\infty}}
& \text{$\sigma=1$.}
\end{cases}
\end{equation*}
For $\sigma=1$ this implies a new identity of the
Rogers--Ramanujan type for the algebra $\D_6$.

\item
If we send $M$ to infinity, replace $q\to 1/q$ and then take the
limit of large  $L$ we obtain for all $k\geq 2$
\begin{multline}\label{X2}
\sum_{r\in\Z^{k-1}_+}
\sideset{}{'}\sum_{m\in\Z_+^7}
\frac{q^{\frac{1}{2}\sum_{a=1}^{k-1}(r_a-r_{a-1})^2}}{(q)_{r_1}}
\Bigl(\prod_{a=2}^{k-1}\qbins{r_{a-1}-r_a+r_{a+1}}{r_a}\Bigr)
q^{\frac{1}{4}mC_{\E_7}m}\qbins{m+n}{m} \\
=\begin{cases}
\chi^{(3k+2,6k-2)}_{(k+1)/2,k}(q) & \text{$k$ odd} \\[3mm]
\chi^{(3k-1,6k+4)}_{k/2,k+1}(q) & \text{$k$ even,}
\end{cases}
\end{multline}
with $r_0=0$, $r_k=r_{k-1}-m_6/2$ and $(m,n)$-system \eqref{mnE7}.
The prime in the sum over $m$ denotes the restriction
that $m_1\equiv m_3\equiv m_5\equiv r_{k-1}\pmod{2}$
and that all other $m_i$ are even.
The identity corresponding to $k=1$ is given by \eqref{E7conj} with $\sigma=0$.
\end{itemize}

\subsection{An ($\A_5$,$\E_6$) series}
Our final conjecture is somewhat more involved than the previous two,
as it is not possible to disentangle the two terms on the left-hand side below
by an appropriate summation restriction on $n$,
\begin{multline}\label{conj3}
\sum_{j\in\Z}\Bigl\{q^{\frac{1}{2}j(48j+2)}\T(L,M,6j,8j)
-q^{\frac{1}{2}(6j+1)(8j+1)}\T(L,M,6j+1,8j+1)\Bigr\} \\
+q^3\sum_{j\in\Z}\Bigl\{q^{\frac{1}{2}j(48j+34)}\T(L,M,6j+2,8j+3)
-q^{\frac{1}{2}(6j+1)(8j+7)}\T(L,M,6j+3,8j+4)\Bigr\} \\
=\sum_{\substack{n\in\Z^5_+ \\ n_1+n_4\equiv n_2+n_5 \pmod{3}\\ 
n_1+n_3+n_5 \equiv L \pmod{2}}}
q^{nC^{-1}_{\A_5}n}\qbin{\frac{1}{2}(L+M+m_3)}{2M}\qbin{m+n}{n},
\end{multline}
with $(m,n)$-system \eqref{mnG} where $\G=\A_5$ (and $p=3$).
Iterating this last conjecture using Theorem~\ref{mainthm} one finds a 
series of $\E_6$-type polynomial identities as follows ($k\geq 1$) 
\begin{multline*}
\sum_{j\in\Z}\Bigl\{q^{\frac{1}{2}j(8(8k+6)j+2)}\T(L,M,(8k+6)j,8j) \\
-q^{\frac{1}{2}(8j+1)((8k+6)j+k+1)}\T(L,M,(8k+6)j+k+1,8j+1)\Bigr\} \\[3mm]
+q^{\frac{3}{2}(3k+2)}\sum_{j\in\Z}\Bigl\{q^{\frac{1}{2}j(8(8k+6)j+48k+34)}
\T(L,M,(8k+6)j+3k+2,8j+3) \\
-q^{\frac{1}{2}(8j+7)((8k+6)j+k+1)}
\T(L,M,(8k+6)j+4k+3,8j+4)\Bigr\} \\
=\sum_{r\in\Z^{k-1}_+}
\Bigl(\prod_{a=0}^{k-2} q^{\frac{1}{2}(r_a-r_{a+1})^2}
\qbins{r_{a-1}-r_a+r_{a+1}}{r_a}\Bigr) \\
\times
\sum_{\substack{n\in\Z^6_+ \\ n_1+n_4\equiv n_2+n_5 \pmod{3}}}
q^{\frac{1}{4}mC_{\E_6}m}\qbins{r_{k-2}-\frac{1}{2}m_6}{r_{k-1}}
\qbins{m+n}{n},
\end{multline*}
where $r_0=L$, $r_{-1}=L+M$ and $m+n=\frac{1}{2}(\I_{\E_6}m+r_{k-1}e_6)$.
Fortunately, the limiting character identities that follow from
this monster are more manageable.
\begin{itemize}
\item
Taking $L+\sigma$ even and letting $M$ and $L$ tend to infinity yields
for positive $k$,
\begin{multline*}
\sum_{n_1,\dots,n_k\geq 0}
\frac{q^{\frac{1}{2}(N_1^2+\cdots+N_k^2)}F^{\D_5}_{n_k;m_{\sigma}}(q)}
{(q)_{n_1}\cdots(q)_{n_{k-1}}(q)_{2n_k}} \\[2mm]
=\begin{cases} \displaystyle
\chi_{\sigma+1,1}^{(3,4)}(q)\chi_{1,(k+1)/2}^{(8,4k+3)}(q)
+\chi_{2-\sigma,1}^{(3,4)}(q)\chi_{7,(k+1)/2}^{(8,4k+3)}(q) 
& \text{$k$ odd} \\[4mm]
\displaystyle
B_{1,k+1;\sigma+k/2}^{(8,8k+6)}(q)+B_{7,k+1;\sigma+k/2+1}^{(8,8k+6)}(q)
& \text{$k$ even,}
\end{cases}
\end{multline*}
with the notation of equation \eqref{Nm}.
The analogous limit taken in \eqref{conj3} leads to
\begin{equation*}
B_{1,1;\sigma}^{(6,8)}(q)+B_{1,7;1-\sigma}^{(6,8)}(q)=
\sum_{\substack{n\in\Z^5_+ \\ n_1+n_4\equiv n_2+n_5 \pmod{3}\\
n_1+n_3+n_5 \equiv \sigma \pmod{2}}}
\frac{q^{nC^{-1}_{\A_5}n}}{(q)_{n}}.
\end{equation*}

\item
If we send $M$ to infinity, replace $q\to 1/q$ and then take the
limit of large  $L$ we obtain for $k\geq 2$,
\begin{multline}\label{X3}
\sum_{r\in\Z^{k-1}_+}
\sideset{}{'}\sum_{m\in\Z_+^6}
\frac{q^{\frac{1}{2}\sum_{a=1}^{k-1}(r_a-r_{a-1})^2}}
{(q)_{r_1}}
\Bigl(\prod_{a=2}^{k-1}
\qbins{r_{a-1}-r_a+r_{a+1}}{r_a}\Bigr)
q^{\frac{1}{4}mC_{\E_6}m}
\qbins{m+n}{m} \\
=\begin{cases}
\chi^{(4k+3,8k-2)}_{(k+1)/2,k}(q)+\chi^{(4k+3,8k-2)}_{(k+1)/2,7k-2}(q)
& \text{$k$ odd} \\[3mm]
\chi^{(4k-1,8k+6)}_{k/2,k+1}(q)+\chi^{(4k-1,8k+6)}_{7k/2-1,k+1}(q)
& \text{$k$ even,}
\end{cases}
\end{multline}
where $m+n=\frac{1}{2}(\I_{\E_6}m+r_{k-1}e_6)$ and $r_0=0$,
$r_k=r_{k-1}-m_6/2$.
The prime in the sum over $m$ denotes the restriction
that $m_1\equiv m_3\equiv m_5\equiv r_{k-1}\pmod{2}$
and that all other $m_i$ are even.
Again $k=1$ is special, corresponding to the $\E_6$ 
conjecture of \cite{KKMM93},
\begin{equation}\label{E6}
\chi^{(6,7)}_{1,1}(q)+\chi^{(6,7)}_{5,1}(q)=
\sum_{\substack{n\in\Z_+^6 \\ n_1+n_4\equiv n_2+n_5 \pmod{3}}}
\frac{q^{nC^{-1}_{\E_6}n}}{(q)_n}.
\end{equation}
\end{itemize}

\section{Discussion}
We hope that the examples presented in the previous section support our
claim that $q$-trinomial coefficients, and their refinement introduced
in this paper are mathematical objects of both depth and elegance.

It is quite intriguing to observe that the $\T$-invariance of 
Theorem~\ref{mainthm} also appears to have physical significance.
In a very recent paper by Dorey, Dunning and Tateo~\cite{DDT00}, new families
of renormalization group flows between $c<1$ conformal field theories 
were proposed. Labelling such a theory by $M(p,p')$,
in accordance with definition~\eqref{Vir} of the Virasoro characters,
Dorey \textit{et al.} conjecture the following flows
\begin{subequations}
\begin{align}\label{a}
M(p,p')+\phi_{21} & \longrightarrow M(p'-p,p') & p<p'<2p & \\
M(p,p')+\phi_{15} & \longrightarrow M(p,4p-p') & 2p<p'<3p & \\
M(p,p')+\phi_{15} & \longrightarrow M(4p-p',p) & 3p<p'<4p &, 
\end{align}
\end{subequations}
where $\phi_{rs}$ is the perturbing operator of the $M(p,p')$ theory.
Using these three flows we find the following chain ending in $M(3,4)$:
$$M(3,4) \stackrel{(\text{c})}{\longleftarrow} M(4,13) 
\stackrel{(\text{a})}{\longleftarrow} M(9,13) 
\stackrel{(\text{b})}{\longleftarrow} M(9,23)
\stackrel{(\text{a})}{\longleftarrow} M(14,23)
\stackrel{(\text{b})}{\longleftarrow} \dots, $$
where (a) denotes equation \eqref{a} etc. But this flow diagram coincides
with the chain of character identities given in \eqref{X} and \eqref{E8}!
In particular, $M(3,4)$ corresponds to the $\E_8$ identity of \eqref{E8},
$M(5n-1,10n+3)$ ($n\geq 1$) corresponds to \eqref{X} for $k=2n$,
and $M(5n+4,10n+3)$ ($n\geq 1$) to \eqref{X} for $k=2n+1$.

In much the same way the flow diagram
$$M(4,5) \stackrel{(\text{c})}{\longleftarrow} M(5,16)
\stackrel{(\text{a})}{\longleftarrow} M(11,16)
\stackrel{(\text{b})}{\longleftarrow} M(11,28)
\stackrel{(\text{a})}{\longleftarrow} M(17,28)
\stackrel{(\text{b})}{\longleftarrow} \dots $$
is in accordance with the chain of character identities given by
\eqref{X2} and \eqref{E7conj}, and
$$M(6,7) \stackrel{(\text{c})}{\longleftarrow} M(7,22)
\stackrel{(\text{a})}{\longleftarrow} M(15,22)
\stackrel{(\text{b})}{\longleftarrow} M(15,38)
\stackrel{(\text{a})}{\longleftarrow} M(23,38)
\stackrel{(\text{b})}{\longleftarrow} \dots $$
is in one-to-one correspondence with \eqref{X3} and \eqref{E6}.

To conclude we mention that more general applications of
our refined $q$-trinomial coefficients will be presented in a 
future paper.
Therein we will also discuss the $\T$-invariance in the broader
context of the Bailey lemma~\cite{Andrews84}, 
trinomial Bailey lemma~\cite{AB98} and (generalized) 
Burge transform~\cite{Burge93,FLW97,SW99}.

\subsection*{Acknowledgements}
This work was supported by a fellowship of the Royal
Netherlands Academy of Arts and Sciences.

\end{document}